\documentclass{amsart}

\usepackage{amsmath} 
\usepackage{amssymb}

\newcommand{\forces}{\Vdash} 
\newcommand{\bV}{{\bf V}} 
\newcommand{\lesdot}{\mathrel{\mathord{<}\!\!\raise 
0.8 pt\hbox{$\scriptstyle\circ$}}} 
\newcommand{\supp}{{\rm supp}}  
\newcommand{\con}{{\mathfrak c}} 
\newcommand{\dominating}{{\mathfrak d\/}} 
\newcommand{\unbounded}{{\mathfrak b}}
 
\newcommand{\baire}{{}^{\textstyle \omega}\omega} 
\newcommand{\iso}{[\omega]^{\textstyle \omega}} 
 
\newcommand{\fseo}{{}^{\textstyle\omega{>}}\omega} 
\newcommand{\conc}{{}^\frown\!}
\newcommand{\lh}{\ell g\/} 
\newcommand{\rest}{{\restriction}}
\newcommand{\dom}{{\rm dom}} 
\newcommand{\rng}{{\rm rng}}
\newcommand{\sqR}{{}^{\textstyle\omega}\mbR}

\newcommand{\cB}{{\mathcal B}}
\newcommand{\cA}{{\mathcal A}}

\newcommand{\cF}{{\mathcal F}}
\newcommand{\cI}{{\mathcal I}}

\newcommand{\bQ}{{\mathbb Q}}  
\newcommand{\dbQ}{{\dot{\mathbb Q}}}

\newcommand{\bP}{{\mathbb P}}
\newcommand{\cP}{{\mathcal P}}
\newcommand{\mbR}{{\mathbb R}}
\newcommand{\cR}{{\mathcal R}}
\newcommand{\cX}{{\mathcal X}}
\newcommand{\cY}{{\mathcal Y}}
\newcommand{\cZ}{{\mathcal Z}}

\newtheorem{theorem}{Theorem}[section] 
\newtheorem{claim}{Claim}[theorem]
 
\newtheorem{proposition}[theorem]{Proposition}

\theoremstyle{definition}
\newtheorem{problem}[theorem]{Problem} 
 
\newtheorem{definition}[theorem]{Definition}

\theoremstyle{remark}
\newtheorem{notation}[theorem]{Notation}
\newtheorem{conclusion}[theorem]{Conclusion}
\newtheorem{remark}[theorem]{Remark}

\begin{document}

\title{The Yellow Cake}

\author{Andrzej Ros{\l}anowski}
\address{Department of Mathematics and Computer Science\\
 Boise State University\\
 Boise ID 83725, USA\\
 and Mathematical Institute of Wroclaw University\\
 50384 Wroclaw, Poland} 
\email{roslanow@math.idbsu.edu}
\urladdr{http://math.idbsu.edu/$\sim$roslanow}
\thanks{The first author thanks the Hebrew University of Jerusalem for
 support during his visit to Jerusalem in Summer'98 when this research was
 done and the KBN (Polish Committee of Scientific Research) for partial
 support through grant 2P03A03114.}  
\author{Saharon Shelah}
\address{Institute of Mathematics\\
 The Hebrew University of Jerusalem\\
 91904 Jerusalem, Israel\\
 and  Department of Mathematics\\
 Rutgers University\\
 New Brunswick, NJ 08854, USA}
\email{shelah@math.huji.ac.il}
\urladdr{http://www.math.rutgers.edu/$\sim$shelah}
\thanks{The research of the second author was partially supported by The
Israel Science Foundation. Publication 686} 

\begin{abstract}
In this paper we consider the following property:
\begin{enumerate}
\item[$(\circledast^{\rm Da})$] For every  function $f:\mbR\times\mbR
\longrightarrow\mbR$ there are functions $g^0_n,g^1_n:\mbR\longrightarrow
\mbR$ (for $n<\omega$) such that 
\[(\forall x,y\in\mbR)(f(x,y)=\sum_{n<\omega}g^0_n(x)g^1_n(y)).\]
\end{enumerate}
We show that, despite some expectation suggested by \cite{Sh:675},
$(\circledast^{\rm Da})$ does not imply ${\bf MA}(\sigma\mbox{--centered})$. 
Next, we introduce cardinal characteristics of the continuum responsible for
the failure of $(\circledast^{\rm Da})$.
\end{abstract}

\maketitle

\section{Introduction}
In the present paper we will consider the following property:
\begin{enumerate}
\item[$(\circledast^{\rm Da})$] For every  function $f:\mbR\times\mbR
\longrightarrow\mbR$ there are functions $g^0_n,g^1_n:\mbR\longrightarrow
\mbR$ (for $n<\omega$) such that 
\[(\forall x,y\in\mbR)(f(x,y)=\sum_{n<\omega}g^0_n(x)g^1_n(y)).\]
\end{enumerate}
Davies \cite{Da74} showed that CH implies $(\circledast^{\rm Da})$ and
Miller \cite[Problem 15.11]{Mixx}, \cite{Mi91} and Ciesielski \cite[Problem
7]{Ci97} asked if $(\circledast^{\rm Da})$ is equivalent to CH. It was shown
in \cite[\S 3]{Sh:675} that the answer is negative. Namely,
\begin{theorem}
\label{675}
\begin{enumerate}
\item (See \cite[3.4]{Sh:675})\quad 
${\bf MA}(\sigma\mbox{-centered})$ implies $(\circledast^{\rm Da})$.
\item (See \cite[3.6]{Sh:675})\quad 
If $\bP$ is the forcing notion for adding $\aleph_2$ Cohen reals then 
$\forces_{\bP}\neg (\circledast^{\rm Da})$.
\end{enumerate}
\end{theorem}
The proof of \cite[Conclusion 3.4]{Sh:675}) strongly used the assumptions
causing an impression that the property $(\circledast^{\rm Da})$ might be
equivalent to ${\bf MA}(\sigma\mbox{-centered})$. 

The first section introduces a strong variant of ccc which is useful in
preserving unbounded families. In the second section we show that
$(\circledast^{\rm Da})$ does not imply ${\bf MA}(\sigma\mbox{-centered})$.
Finally, the in next section we show the combinatorial heart of
\cite[Proposition 3.6]{Sh:675} and we introduce cardinal characteristics of
the continuum closely related to the failure of $(\circledast^{\rm Da})$.
\medskip

\noindent{\bf Notation}\qquad Most of our notation is standard and
compatible with that of classical textbooks on Set Theory (like
Bartoszy\'nski Judah \cite{BaJu95}). However in forcing we keep the
convention that {\em a stronger condition is the larger one}.

\begin{notation}
\label{notacja}
\begin{enumerate}
\item For two sequences $\eta,\nu$ we write $\nu\vartriangleleft\eta$ whenever
$\nu$ is a proper initial segment of $\eta$, and $\nu\trianglelefteq\eta$ when
either $\nu\vartriangleleft\eta$ or $\nu=\eta$. The length of a sequence
$\eta$ is denoted by $\lh(\eta)$.
\item The set of rationals is denoted by $\bQ$ and the set of reals is
called $\mbR$. The cardinality of $\mbR$ is called $\con$ (and it is refered
to as the the continuum). The dominating number (the minimal size of a
dominating family in $\baire$ in the ordering of eventual dominance) is
denoted by $\dominating$ and the unbounded number (the minimal size of an
unbounded family in that order) is called $\unbounded$.
\item The quantifiers $(\forall^\infty n)$ and $(\exists^\infty n)$ are
abbreviations for  
\[(\exists m\in\omega)(\forall n>m)\quad\mbox{ and }\quad(\forall m\in\omega)
(\exists n>m),\] 
respectively.
\item For a forcing notion $\bP$, $\Gamma_\bP$ stands for the canonical
$\bP$--name for the generic filter in $\bP$. With this one exception, all
$\bP$--names for objects in the extension via $\bP$ will be denoted with a dot
above (e.g.~$\dot{A}$, $\dot{f}$).
\end{enumerate}
\end{notation}

\section{$\cF$--sweet forcing notion}

\begin{definition}
\label{spread}
An uncountable family $\cF\subseteq\baire$ is {\rm spread} if
\begin{enumerate}
\item[$(\boxtimes)$] for each $k^*,n^*<\omega$ and a sequence $\langle
f_{\alpha,n}:\alpha<\omega_1,\ n<n^*\rangle$ of pairwise distinct elements
of $\cF$ there are an increasing sequence $\langle\alpha_i:i<\omega\rangle
\subseteq\omega_1$ and an integer $k>k^*$ such that
\[(\forall i<\omega)(\forall n<n^*)(f_{\alpha_i,n}(k)<f_{\alpha_{i+1},n}
(k)).\]
\end{enumerate}
\end{definition}

\begin{remark}
\begin{enumerate}
\item Note that if an uncountable family $\cF\subseteq\baire$ has the
property that its every uncountable subfamily is unbounded on every $K\in
\iso$ then $\cF$ is spread.
\item If $\kappa$ is uncountable and one adds $\kappa$ many Cohen reals
$\langle c_\alpha:\alpha<\kappa\rangle\subseteq\baire$ then $\{c_\alpha:
\alpha<\kappa\}$ is a spread family. 
\item If there is a spread family then $\unbounded=\aleph_1$ (so in
particular ${\bf MA}_{\aleph_2}(\sigma\mbox{--centered})$ fails).
\end{enumerate}
\end{remark}

\begin{definition}
\label{sweet}
Let $\cF\subseteq\baire$ be a spread family. A forcing notion $\bP$ is {\em
$\cF$--sweet} if the following condition is satisfied:
\begin{enumerate}
\item[$(\boxplus)^{\cF}_{\rm sweet}$] for each sequence $\langle p_\alpha:
\alpha<\omega_1\rangle\subseteq\bP$ there are $A\in [\omega]^{\textstyle
\aleph_1}$, $k^*<\omega$ and a sequence $\langle f_{\alpha,n}:n<n^*,\
\alpha\in A \rangle\subseteq\cF$ such that $(\alpha,n)\neq (\alpha',n')\
\Rightarrow\ f_{\alpha,n}\neq f_{\alpha',n'}$ and 
\begin{enumerate}
\item[$(\oplus)$] if $\langle\alpha_i:i<\omega\rangle$ is an increasing
sequence of elements of $A$ such that for some $k\in (k^*,\omega)$ 
\[(\forall i<\omega)(\forall n< n^*)(f_{\alpha_i,n}(k)<f_{\alpha_{i+1},n}
(k))\]
then there is $p\in\bP$ such that $p\forces(\exists^\infty i\in\omega)(
p_{\alpha_i}\in\Gamma_{\bP})$. 
\end{enumerate}
\end{enumerate}
\end{definition}

\begin{proposition}
\label{swspr}
Assume that $\cF\subseteq\baire$ is a spread family and $\bP$ is an
$\cF$--sweet forcing notion. Then 
\[\forces_{\bP}\mbox{`` $\cF$ is a spread family ''.}\]
\end{proposition}

\begin{proof}
First note that easily $\cF$--sweetness implies the ccc.

Suppose that $k^+<\omega$, $\langle\dot{f}_{\alpha,n}:\alpha<\omega_1,\ n<
n^+\rangle$ are $\bP$--names for elements of $\cF$, $p\in\bP$ and 
\[p\forces_{\bP}(\forall\alpha,\alpha'<\omega_1)(\forall n,n'<n^+)((\alpha,n)
\neq (\alpha',n')\ \Rightarrow\ \dot{f}_{\alpha,n}\neq\dot{f}_{\alpha',n'}).\]
For $\alpha<\omega_1$ choose conditions $p_\alpha\geq p$ and functions
$f_{\alpha,n}\in\cF$ (for $n<n^+$) such that $p_\alpha\forces (\forall
n<n^+)(\dot{f}_{\alpha,n}=f_{\alpha,n})$. Passing to a subsequence, we
may assume that  
\[(\alpha,n)\neq (\alpha',n')\ \Rightarrow\ f_{\alpha,n}\neq
f_{\alpha',n'}.\]
Choose $k^*>k^+$, a set $A\in [\omega]^{\textstyle\aleph_1}$ and a
sequence $\langle f_{\alpha,n}: \alpha\in A,\ n^+\leq n<n^*\rangle$ as
guaranteed by $(\boxplus)^\cF_{\rm sweet}$ of \ref{sweet} for $\langle
p_\alpha:\alpha<\omega_1\rangle$ (note that here, for notational
convenience, we use the interval $[n^+,n^*)$ instead of $n^*$ there).
Shrinking the set $A$ and possibly decreasing $n^*$ (and reenumerating
$f_{\alpha,n}$'s) we may assume that all functions in appearing in $\langle
f_{\alpha,n}: \alpha\in A,\ n<n^*\rangle$ are distinct. By $(\boxtimes)$ of
\ref{spread} we find $k>k^*$ and an increasing sequence $\langle\alpha_i:
i<\omega\rangle\subseteq A$ such that 
\[(\forall i<\omega)(\forall n<n^*)(f_{\alpha_i,n}(k)<f_{\alpha_{i+1},n}
(k)).\]
But it follows from $(\oplus)$ of \ref{sweet} that now we can find a
condition $q\in\bP$ such that $q\forces(\exists^\infty i\in\omega)(
p_{\alpha_i}\in\Gamma_{\bP})$. As all conditions $p_\alpha$ are stronger
than $p$ we may demand that $q\geq p$. Now use the choice of the
$p_{\alpha_i}$'s and $f_{\alpha_i,n}$ (for $n<n^+$) to finish the proof. 
\end{proof}

\begin{theorem}
\label{pressweet}
Assume $\cF$ is a spread family. Let $\langle\bP_\alpha,\dbQ_\alpha:\alpha<
\gamma\rangle$ be a finite support iteration of forcing notions such that
for each $\alpha<\gamma$ we have
\begin{enumerate}
\item $\forces_{\bP_\alpha}$`` $\cF$ is spread '', and
\item $\forces_{\bP_\alpha}$`` $\dbQ_\alpha$ is $\cF$--sweet ''.
\end{enumerate}
Then $\bP_\gamma$ is $\cF$--sweet (and consequently, $\forces_{
\bP_\gamma}$`` $\cF$ is a spread family ''). 
\end{theorem}

\begin{proof}
We show this by induction on $\gamma$.
\medskip

\noindent{\sc Case 1:}\qquad $\gamma=\beta+1$\\
Let $\langle p_\alpha:\alpha<\omega_1\rangle\subseteq\bP_{\beta+1}$. Take a
condition $p^*\in\bP_\beta$ such that 
\[p^*\forces_{\bP_\beta}\mbox{`` }\{\alpha<\omega_1:p_\alpha\rest\beta\in
\Gamma_{\bP_\beta}\}\mbox{ is uncountable ''}\]
(there is one by the ccc). Next, use the assumption that $\dbQ_\beta$ is
$\cF$--sweet and get $\bP_\beta$--names $\dot{A}\in [\omega_1]^{\textstyle
\aleph_1}$ and $\dot{k}^*$, $\dot{n}^*$ and $\langle \dot{f}_{\alpha,n}:
\alpha\in\dot{A},\ n<\dot{n}^*\rangle\subseteq\cF$ such that the condition
$p^*$ forces that they are as guaranteed by $(\boxplus)^\cF_{\rm sweet}$ of
\ref{sweet} for the sequence $\langle p_\alpha(\beta):\alpha<\omega_1,\
p_\alpha\rest\beta\in\Gamma_{\bP_\beta}\rangle$. 

Let $A'$ be the set of all $\alpha<\omega_1$ such that there is a condition
stronger than both $p^*$ and $p_\alpha\rest\beta$ which forces that
$p_\alpha(\beta)$ is in $\dot{A}$. Clearly $|A'|=\aleph_1$. For each
$\alpha\in A'$ choose a condition $q_\alpha\in\bP_\beta$ stronger than both
$p^*$ and $p_\alpha\rest\beta$ which forces that $p_\alpha(\beta)\in\dot{A}$
and decides the values of $\dot{k}^*$, $\dot{n}^*$ and $\langle\dot{f}_{
\alpha,n}: n<\dot{n}^*\rangle$. Next we may choose $A''\in [A']^{\textstyle
\aleph_1}$, $k^*,n^*$ and $\langle f_{\alpha,n}:\alpha\in A'',\ n<n^*\rangle
\subseteq\cF$ such that (for each $\alpha\in A''$ and $n<n^*$) $q_\alpha
\forces$`` $\dot{k}^*=k^*\ \&\ \dot{n}^*=n^*\ \&\ \dot{f}_{\alpha,n}=
f_{\alpha,n}$ ''. Moreover we may demand that the $f_{\alpha,n}$'s are
pairwise distinct (for $\alpha\in A''$, $n<n^*$).

Apply the inductive hypothesis to the sequence $\langle q_\alpha:\alpha\in
A''\rangle$ (and $\bP_\beta$) to get $A\in [A'']^{\textstyle\aleph_1}$,
$k^+$, $n^+>n^*$ and $\langle f_{\alpha,n}:\alpha\in A,\ n^*\leq n<n^+
\rangle$. For simplicity we may assume that there are no repetitions in the
sequence $\langle f_{\alpha,n}:\alpha\in A,\ n<n^*\rangle$ (we may shrink
$A$ and decrease $n^*$ reenumerating $f_{\alpha,n}$'s suitably). We claim
that this sequence and $\max\{k^*,k^+\}$ satisfy the demand in $(\oplus)$ if
\ref{sweet}. So suppose that $\langle\alpha_i:i<\omega\rangle$ is an
increasing sequence of elements of $A$ such that for some $k>k^*,k^+$ we
have  
\[(\forall i<\omega)(\forall n<n^+)(f_{\alpha_i,n}(k)<f_{\alpha_{i+1},n}
(k)).\]
Clearly, by our choices, we find a condition $p^+\in\bP_\beta$ stronger than
$p^*$ such that $p^+\forces(\exists^\infty i\in\omega)(q_{\alpha_i}\in
\Gamma_{\bP_\beta})$. Next, in $\bV^{\bP_\beta}$, we look at the sequence
$\langle p_{\alpha_i}(\beta):q_{\alpha_i}\in\Gamma_{\bP_\beta},\ i<\omega
\rangle$. We may find a $\bP_\beta$--name $p^+(\beta)$ such that ($p^+$
forces that)  
\[p^+(\beta)\forces_{\dbQ_\beta}(\exists^\infty i\in\omega)(q_{\alpha_i}\in
\Gamma_{\bP_\beta}\ \&\ p_{\alpha_i}(\beta)\in\Gamma_{\dbQ_\beta}).\]
Look at the condition $p^+\conc p^+(\beta)$. 
\medskip

\noindent{\sc Case 2:}\qquad $\gamma$ is a limit ordinal.\\
If $\langle p_\alpha:\alpha<\omega_1\rangle\subseteq\bP_\gamma$ then, under
the assumption of the current case, for some $A\in [\omega_1]^{\textstyle
\aleph_1}$ and $\delta<\gamma$, the sets $\{\supp(p_\alpha)\setminus\delta:
\alpha\in A\}$ are pairwise disjoint. Apply the inductive hypothesis to
$\bP_\delta$ and the sequence $\langle p_\alpha\rest\delta:\alpha\in A
\rangle$. 
\end{proof}

\begin{conclusion}
\label{getMA}
Suppose that $\kappa>\aleph_1$ is a regular cardinal such that $\kappa^{{<}
\kappa}=\kappa$ and $(\forall\mu<\kappa)(\mu^{\aleph_0}<\kappa)$. Then there
is a ccc forcing notion $\bP$ of size $\kappa$ such that 
\[\forces_{\bP}\mbox{`` there is a spread family $\cF\subseteq\baire$ of
size $\kappa$ } \&\ \con=\kappa\ \&\ {\bf MA}(\cF\mbox{--sweet})\mbox{ ''.}\]
\end{conclusion}

\begin{proof}
First note that if $\bP$ is an $\cF$--sweet forcing notion,
$\cI_\xi\subseteq\bP$ (for $\xi<\mu<\kappa$) are dense subsets of $\bP$ and
$p\in\bP$ then, under our assumptions, there is a set $\bP^*\subseteq\bP$ of
size less than $\kappa$ such that $p\in\bP^*$ and
\begin{itemize}
\item if $p,q\in\bP^*$ are incompatible in $\bP^*$ then they are
incompatible in $\bP$,
\item if $\langle p_i:i<\omega\rangle\subseteq\bP^*$ is not a maximal
antichain in $\bP$ then it is not in $\bP^*$, 
\item for each $\xi<\mu$ the intersection $\cI_\xi\cap\bP^*$ is dense in
$\bP^*$. 
\end{itemize}
(Thus $\bP^*\lesdot\bP$ and so it is $\cF$--sweet.)

Now, using standard bookkeeping arguments, build a finite support iteration
$\langle \bP_\alpha,\dbQ_\alpha:\alpha<\kappa\rangle$ such that 
\begin{enumerate}
\item $\bQ_0$ is the forcing notion adding $\kappa$ many Cohen real $\cF=
\langle f_\alpha:\alpha<\kappa\rangle\subseteq\baire$ (with finite
conditions), [so in $\bV^{\bQ_0}$, the family $\cF$ is spread]
\item for each $\alpha<\kappa$, $\forces_{\bP_{1+\alpha}}$``
$\dbQ_{1+\alpha}$ is a $\cF$--sweet forcing notion of size $<\kappa$ '',
\item if $\dbQ$ is a $\bP_\kappa$--name for a $\cF$--sweet forcing notion of
size $<\kappa$ then for $\kappa$ many $\alpha<\kappa$, $\dbQ$ is a
$\bP_\alpha$--name and $\forces_{\bP_\alpha}\dbQ=\dbQ_\alpha$.
\end{enumerate}
It follows from \ref{pressweet} that in $\bV^{\bP_\alpha}$ (for $0<\alpha
\leq\kappa$) the family $\cF$ is spread, so there are no problems with
carrying out the construction. Easily $\bP_\kappa$ is as required. 
\end{proof}

\begin{remark}
Note the similarity of ${\bf MA}(\cF\mbox{--sweet})$ to the methods used in
\cite[\S 4]{Sh:98}.
\end{remark}

\section{More on Davies' Problem}
The aim of this section is to show that $(\circledast^{\rm Da})$ does not
imply ${\bf MA}(\sigma\mbox{--centered})$. 

Let $\langle\nu_n:n<\omega\rangle$ be an enumeration of $\fseo$ such that
$\lh(\nu_n)\leq n$. For distinct $\rho_0,\rho_1\in\baire$ let
$\delta(\rho_0,\rho_1)=1+\max\{m:\nu_m\vartriangleleft\rho_0\ \&\ \nu_m
\vartriangleleft\rho_1\}$. (Note that $\rho_0\rest\delta(\rho_0,\rho_1)\neq
\rho_1\rest\delta(\rho_0,\rho_1)$.)

Assume that there exists a spread family of size $\con$ and let $\cF=\langle
\rho_\alpha:\alpha<\con\rangle\subseteq\baire$ be such a family (later we
will choose the one coming from adding $\kappa$ many Cohen reals).

\begin{definition} 
\label{aproksymacje} 
Let $\zeta<\con$ be an ordinal and let $f:\zeta\times\zeta\longrightarrow
\mbR$.
\begin{enumerate}
\item {\em A $\zeta$--approximation} is a sequence $\bar{g}=\langle
g^\ell_\eta:\ell<2,\ \eta\in\fseo\rangle$ such that:
\begin{enumerate}
\item[(a)] $g^\ell_\eta:\zeta\longrightarrow\bQ$ (for $\ell<2$, $\eta\in
\fseo$),  
\item[(b)] if $\alpha<\zeta$ then $(\forall\beta<\con)(\exists^\infty k\in
\omega)(g^\ell_{\rho_\alpha\rest k}(\alpha)\neq 0\ \&\ \nu_k
\vartriangleleft\rho_\beta)$,
\item[(c)] if $\alpha<\zeta$, $\eta\in\fseo$ and neither $\eta$ nor $\nu_{
\lh(\eta)}$ is an initial segment of $\rho_\alpha$, then $g^0_\eta(\alpha)
=g^1_\eta(\alpha)=0$.
\end{enumerate}
\item If $\zeta_0<\zeta_1$ and $\bar{g}^k=\langle g^{\ell,k}_\eta:\ell<2,\
\eta\in\fseo\rangle$ (for $k=0,1$) are $\zeta_k$--approximations such that
$g^{\ell,0}_\eta\subseteq g^{\ell,1}_\eta$ (for all $\ell<2$ and $\eta\in
\fseo$) then we say that {\em $\bar{g}^1$ extends $\bar{g}^0$} (in short:
$\bar{g}^0\preceq\bar{g}^1$).  
\item We say that a $\zeta$--approximation $\bar{g}$ {\em agrees with the
function $f$} if   
\[(\forall\alpha,\beta<\zeta)\big(f(\alpha,\beta)=\sum_{\eta\in{}^{\omega{>}}
\omega}g^0_\eta(\alpha)\cdot g^1_\eta(\beta)\mbox{ and the series converges
absolutely\/}\big).\]  
\end{enumerate}
\end{definition}

\begin{proposition}
\label{clear}
If $\bar{g}^\xi$ are $\zeta_\xi$--approximations (for $\xi<\xi^*$) such that
the sequence $\langle\bar{g}^\xi:\xi<\xi^*\rangle$ is $\preceq$--increasing
and $\zeta_{\xi^*}=\bigcup\limits_{\xi<\xi^*} \zeta_\xi$ then there is a
$\zeta_{\xi^*}$--approximation $\bar{g}^{\xi^*}$ such that
$(\forall\xi<\xi^*)(\bar{g}^\xi\preceq\bar{g}^{\xi^*})$. Moreover, if
$f:\zeta_{\xi^*}\times\zeta_{\xi^*}\longrightarrow\mbR$ and each
$\bar{g}^\xi$ agrees with $f\rest(\zeta_\xi\times\zeta_\xi)$ then
$\bar{g}^{\xi^*}$ agrees with $f$. 
\end{proposition}

Thus if we want to show that $(\circledast^{\rm Da})$ holds we may take a
function $f:\con\times\con\longrightarrow\mbR$ (it should be clear that we
may look at functions of that type only) and try to build a
$\preceq$--increasing sequence $\langle\bar{g}^\xi:\xi<\con\rangle$ of
approximations. If we make sure that $\bar{g}^\xi$ is a $\xi$--approximation
that agrees with $f\rest(\xi\times\xi)$ then the limit $\bar{g}^{\con}$ of
$\bar{g}^\xi$'s will give us witnesses for $f$. (Note that by the absolute
convergence demand in \ref{aproksymacje}(3) we do not have to worry about
the order in the series.) At limit stages of the construction we use
\ref{clear}, but problems may occur at some successor stage. Here we need to
use forcing. 

\begin{definition}
\label{forcing}
Assume that $\zeta<\con$ is an ordinal, and $f:(\zeta+1)\times(\zeta+1)
\longrightarrow\mbR$. Let $\bar{g}=\langle g^\ell_\eta:\ell<2,\ \eta\in
\fseo\rangle$ be a $\zeta$--approximation which agrees with $f\rest\zeta
\times\zeta$. We define a forcing notion $\bP^{\bar{g},\zeta}_{f}$ as
follows: 
\smallskip

\noindent{\bf a condition} is a tuple $p=\langle Z^p,j^p,\langle r^p_{\ell,
\eta}:\ell<2,\ \eta\in{}^{\textstyle j^p{>}}\omega\rangle\rangle$ such that 
\begin{enumerate}
\item[$(\alpha)$] $j^p<\omega$ and $Z^p$ is a finite subset of $\zeta$,
$r^p_{\ell,\eta}\in\bQ$ (for $\ell<2$, $\eta\in {}^{\textstyle
j^p{>}}\omega$),  
\item[$(\beta)$]  the set $\{\eta\in {}^{\textstyle j^p{>}}\omega:r^p_{0,
\eta}\neq 0\mbox{ or } r^p_{1,\eta}\neq 0\}$ is finite, and if $\eta\in
{}^{\textstyle j^p{>}}\omega$ and neither $\eta$ nor $\nu_{\lh(\eta)}$ is an
initial segment of $\rho_\zeta$ then $r^p_{\ell,\eta}=0$,
\item[$(\gamma)$] if $\alpha\in Z^p$ then 
\[\begin{array}{l}
|f(\alpha,\zeta)-\sum\{g^0_\eta(\alpha)\cdot r^p_{1,\eta}:\eta\in {}^{
\textstyle j^p{>}}\omega\}|<2^{-j^p},\\
|f(\zeta,\alpha)-\sum\{r^p_{0,\eta}\cdot g^1_\eta(\alpha):\eta\in {}^{
\textstyle j^p{>}}\omega\}|<2^{-j^p},\quad\mbox{ and}\\
|f(\zeta,\zeta)-\sum\{r^p_{0,\eta}\cdot r^p_{1,\eta}:\eta\in {}^{
\textstyle j^p{>}}\omega\}|<2^{-j^p}
  \end{array}\]
(note that by demand $(\beta)$ all the sums above are finite),
\item[$(\delta)$] if $\alpha,\beta\in Z^p\cup\{\zeta\}$ are distinct then
$\delta(\rho_\alpha,\rho_\beta)<j^p$;  
\end{enumerate}
{\bf the order} is defined by $p\leq q$ if and only if
\begin{enumerate}
\item[(a)] $j^p\leq j^q$, $Z^p\subseteq Z^q$ and $r^p_{\ell,\eta}=r^q_{
\ell,\eta}$ for $\eta\in {}^{\textstyle j^p{>}}\omega$, $\ell<2$,
\item[(b)] if $\alpha\in Z^p$ then
\[\begin{array}{l}
\displaystyle\sum\{|r^p_{0,\eta}\cdot g^1_\eta(\alpha)|:\eta\in {}^{
\textstyle j^q{>}}\omega\setminus{}^{\textstyle j^p{>}}\omega\}<4\frac{1-
2^{j^p-j^q}}{2^{j^p-1}},\\
\displaystyle\sum\{|g^0_\eta(\alpha)\cdot r^p_{1,\eta}|:\eta\in {}^{
\textstyle j^q{>}}\omega\setminus{}^{\textstyle j^p{>}}\omega\}<4\frac{1-
2^{j^p-j^q}}{2^{j^p-1}},\quad\mbox{ and}\\
\displaystyle\sum\{|r^p_{0,\eta}\cdot r^p_{1,\eta}|:\eta\in {}^{\textstyle
j^q{>}}\omega\setminus{}^{\textstyle j^p{>}}\omega\}<4\frac{1-2^{j^p-j^q}}{
2^{j^p-1}}.\\
  \end{array}\]
\end{enumerate}
\end{definition}

\begin{proposition}
\label{forissweet}
Suppose that $\zeta<\con$, $f:(\zeta+1)\times(\zeta+1)\longrightarrow\mbR$
and $\bar{g}$ is a $\zeta$--approximation that agrees with $f\rest\zeta
\times\zeta$. Then:
\begin{enumerate}
\item $\bP^{\bar{g},\zeta}_{f}$ is a (non-trivial) $\cF$--sweet forcing
notion of size $|\zeta|+\aleph_0$.
\item In $\bV^{\bP^{\bar{g},\zeta}_{f}}$, there is a
$(\zeta+1)$--approximation $\bar{g}^*$ such that $\bar{g}\prec\bar{g}^*$ and
$\bar{g}^*$ agrees with $f$. 
\end{enumerate}
\end{proposition}

\begin{proof}
(1)\qquad First note that $(\bP^{\bar{g},\zeta}_f,\leq)$ is a partial
order and easily $\bP^{\bar{g},\zeta}_f\neq\emptyset$ (remember that $Z^p$
may be empty). Before we continue let us show the following claim that will
be used later too.

\begin{claim}
\label{cl1}
For each $j<\omega$, $\xi<\zeta$ and $\rho\in\baire$ the sets 
\[\begin{array}{lcl}
\cI^j&\stackrel{\rm def}{=}&\{p\in\bP^{\bar{g},\zeta}_f: j^p\geq j\},\\
\cI_\xi&\stackrel{\rm def}{=}&\{p\in\bP^{\bar{g},\zeta}_f: \xi\in Z^p\},
\quad\mbox{ and}\\
\cI^j_\rho&\stackrel{\rm def}{=}&\{p\in\bP^{\bar{g},\zeta}_f: j<j^p\ \&\
(\forall\ell<2)(\exists k\in (j,j^p))(r^p_{\ell,\rho_\zeta\rest j}\neq 0\
\&\ \nu_j\vartriangleleft\rho)\}\\
  \end{array}\]
are dense subsets of $\bP^{\bar{g},\zeta}_f$. 
\end{claim}

\begin{proof}[Proof of the claim]
Let $j<\omega$, $\xi<\zeta$, $\rho\in\baire$ and $p\in\bP^{\bar{g},
\zeta}_f$. 

If $j\leq j^p$ then $p\in\cI^j$, so suppose that $j^p<j$. Let $\langle\xi_m:
m<m^*\rangle$ enumerate $Z^p$. Choose pairwise distinct $\langle j_{\ell,m}:
\ell<2, m<m^*\rangle\subseteq (j,\omega)$ such that $\nu_{j_{\ell,m}}
\vartriangleleft \rho_\zeta$ and $g^\ell_{\rho_{\xi_m}\rest j_{\ell,m}}(
\xi_m)\neq 0$ (remember \ref{aproksymacje}(1b)). Fix $j^*>j$ such that
$\nu_{j^*}$ is not an initial segment of any $\rho_{\xi_m}$ (for $m<m^*$). 
Let $j^q=j+\max\{j_{\ell,m}:\ell<2,\ m<m^*\}+j^*$, $Z^q=Z^p$  and define
$r^q_{0,\eta},r^q_{1,\eta}$ as follows.
\begin{enumerate}
\item If $\eta\in {}^{\textstyle j^p{>}}\omega$ then $r^q_{\ell,\eta}=r^p_{
\ell,\eta}$.
\item If $\eta\in {}^{\textstyle j^q{>}}\omega\setminus {}^{\textstyle
j^p{>}}\omega\setminus \{\rho_{\xi_m}\rest j_{\ell,m}: m<m^*\}\setminus
\{\rho_\zeta\rest j^*\}$, $\ell<2$ then $r^q_{1-\ell,\eta}=0$.
\item If $\eta=\rho_\zeta\rest j^*$ then $r^q_{0,\eta},r^q_{1,\eta}\in\bQ
\setminus\{0\}$ are such that $|r^q_{0,\eta}\cdot r^q_{1,\eta}|<2^{-j^p}$
and 
\[|f(\zeta,\zeta)-\sum\{r^p_{0,\nu}\cdot r^p_{1,\nu}:\nu\in {}^{\textstyle
j^p{>}}\omega\}- r^q_{0,\eta}\cdot r^q_{1,\eta}|<2^{-2j^q}.\]
\item If $\eta=\rho_{\xi_m}\rest j_{0,m}$, $m<m^*$ then $r^q_{1,\eta}\in\bQ$
is such that $|g^0_\eta(\xi_m)\cdot r^q_{1,\eta}|<2^{-j^p}$ and 
\[|f(\xi_m,\zeta)-\sum\{g^0_\nu(\xi_m)\cdot r^p_{1,\nu}:\nu\in{}^{\textstyle
j^p{>}}\omega\}-g^0_\eta(\xi_m)\cdot r^q_{1,\eta}|<2^{-2j^q};\]
if $\eta=\rho_{\xi_m}\rest j_{1,m}$, $m<m^*$ then $r^q_{0,\eta}\in\bQ$ is
such that $|r^q_{0,\eta}\cdot g^1_\eta(\xi_m)|<2^{-j^p}$ and 
\[|f(\zeta,\xi_m)-\sum\{r^p_{0,\nu}\cdot g^1_\nu(\xi_m):\nu\in{}^{\textstyle
j^p{>}}\omega\}-r^q_{0,\eta}\cdot g^1_\eta(\xi_m)|<2^{-2j^q}.\]
\end{enumerate}
One easily checks that $q=\langle Z^q,j^q,\langle r^q_{\ell,\eta}:\ell<2,\
\eta\in{}^{\textstyle j^q{>}}\omega\rangle\rangle$ is a condition in
$\bP^{\bar{g},\zeta}_f$ stronger than $p$ (and $q\in\cI^j$).

Now suppose that $\xi\notin Z^p$. Take $j_0>j^p$ such that $(\forall\alpha
\in Z^p\cup\{\zeta\})(\delta(\xi,\alpha)<j_0)$. Let $\langle\xi_m:m<m^*
\rangle$ enumerate $Z^p\cup\{\xi\}$ and let $\langle j_{\ell,m}:\ell<2,\
m<m^*\rangle\subseteq (j_0,\omega)$ be pairwise distinct and such that
$\nu_{j_{\ell,m}}\vartriangleleft\rho_\zeta\ \&\ g^\ell_{\xi_m\rest
j_{\ell,m}}(\xi_m)\neq 0$. Let $j^*>j^p$ be such that $\nu_{j^*}$ is not an
initial segment of any $\rho_{\xi_m}$. Put $Z^q=Z^p\cup\{\xi\}$, $j^q=j^p+
\max\{j_{\ell,m}:\ell<2,\ m<m^*\}+j^*$, and define $r^q_{\ell,\eta}$ like
before, with one modification. If $\xi_m=\xi$ and $\eta=\rho_\xi\rest j_{0,
m}$ then $r^q_{1,\eta}\in\bQ$ is such that $|f(\xi,\zeta)-g^0_\eta(\xi)\cdot
r^q_{1,\eta}|<2^{-2j^q}$; if $\xi_m=\xi$ and $\eta=\rho_\xi\rest j_{1,m}$
then $r^q_{0,\eta}\in\bQ$ is such that $|f(\zeta,\xi)-r^q_{0,\eta}\cdot
g^1_\eta(\xi)|<2^{-2j^q}$. 

Similarly one builds a condition $q\in\cI^j_\rho$ stronger than $p$
(just choose $j^*$ suitably).
\end{proof}
Now we are going to show that $\bP^{\bar{g},\zeta}_f$ is $\cF$--sweet. So
suppose that $\langle p_\alpha:\alpha<\omega\rangle\subseteq\bP^{\bar{g},
\zeta}_f$. Choose $A\in[\omega_1]^{\textstyle\aleph_1}$ such that 
\begin{itemize}
\item $\langle Z^{p_\alpha}:\alpha\in A\rangle$ forms a $\Delta$--system
with kernel $Z$,
\item for each $\alpha,\beta\in A$, $|Z^{p_\alpha}|=|Z^{p_\beta}|$,
$j^{p_\alpha}=j^{p_\beta}$ and 
\[\langle r^{p_\alpha}_{\ell,\eta}:\ell<2,\ \eta\in{}^{\textstyle
j^{p_\alpha}{>}}\omega\rangle=\langle r^{p_\beta}_{\ell,\eta}:\ell<2,\
\eta\in{}^{\textstyle j^{p_\beta}{>}}\omega\rangle\]
(remember \ref{forcing}($\beta$)),
\item if $\alpha,\beta\in A$ and $\pi:Z^{p_\alpha}\longrightarrow Z^{
p_\beta}$ is the order preserving bijection then $\pi\rest Z$ is the
identity on $Z$ and $(\forall\xi\in Z^{p_\alpha})(\rho_\xi\rest
j^{p_\alpha}=\rho_{\pi(\xi)}\rest j^{p_\beta})$.
\end{itemize}
Let $k^*=j^{p_\alpha}$, $n^*=|Z^{p_\alpha}\setminus Z|$ for some
(equivalently: all) $\alpha\in A$. For $\alpha\in A$ let $\langle f_{\alpha,
n}:n<n^*\rangle$ enumerate $\{\rho_\xi:\xi\in Z^{p_\alpha}\setminus Z\}$. 
Clearly there are no repetitions in $\langle f_{\alpha,n}:n<n^*,\ \alpha\in
A\rangle$. We claim that this sequence is as required in $(\oplus)$ of
\ref{sweet}. So suppose that $\langle\alpha_i:i<\omega\rangle\subseteq A$ is
an increasing sequence such that for some $k>k^*$ we have
\[(\forall i<\omega)(\forall n<n^*)(f_{\alpha_i,n}(k)<f_{\alpha_{i+1},n}
(k)).\]
Passing to a subsequence we may additionally demand that for each $m<k$, for
every $n<n^*$, the sequence $\langle f_{\alpha_i,n}(m): i<\omega\rangle$ is
either constant or strictly increasing. For $n<n^*$ let $k_n\geq j^p$ be
such that the sequence $\langle f_{\alpha_i,n}\rest k_n:i<\omega\rangle$ is
constant but the sequence $\langle f_{\alpha_i,n}(k_n):i<\omega\rangle$ is
strictly increasing. Take $j>k$ such that if $\nu_m\trianglelefteq
f_{\alpha_i,n}\rest k_n$, $n<n^*$ then $m<j$. Fix an enumeration $\langle
\xi_m:m<m^*\rangle$ of $Z^{p_{\alpha_0}}$ (so $m^*=|Z|+n^*$) and choose
$j^*,j_{\ell,m}>j+2$ with the properties as in the first part of the proof
of \ref{cl1} (with $p_{\alpha_0}$ in the place of $p$ there). Put $Z^q=Z^{
p_{\alpha_0}}$ and define $j^q,r^q_{\ell,\eta}$ exactly as there (so, in
particular, for each $\eta\in{}^{\textstyle j{>}}\omega\setminus
{}^{\textstyle j_{p_{\alpha_0}}{>}}\omega$ we have $r^q_{\ell,\eta}=0$). We
claim that $q\forces(\exists^\infty i\in\omega)(p_{\alpha_i}\in\Gamma_{
\bP^{\bar{g},\zeta}_f})$. So suppose that $q'\geq q$, $i_0<\omega$. Choose
$i>i_0$ such that for each $n<n^*$ and $k'>k_n$, if $\nu_m=f_{\alpha_i,n}
\rest k'$ then $m>j^{q'}$. Moreover, we demand that if $k_n<k'<j^{q'}$,
$n<n^*$ then $r^{q'}_{0,f_{\alpha_i,n}\rest k'}=r^{q'}_{1,f_{\alpha_i,n}
\rest k'}=0$ (remember \ref{forcing}($\beta$)). Then we have the effect that
\[(\forall\eta\in{}^{\textstyle j^{q'}{>}}\omega\setminus{}^{\textstyle
j^{p_{\alpha_i}{>}}}\omega)(\forall\ell<2)(\forall\xi\in Z^{p_{\alpha_i}}
\setminus Z)(r^{q'}_{\ell,\eta}\cdot g^1_\eta(\xi)=g^0_\eta(\xi)\cdot
r^{q'}_{1,\eta}=0)).\]
So we may proceed as in the proof of \ref{cl1} and build a condition $q^+$
stronger than both $q'$ and $p_{\alpha_i}$. 
\medskip

\noindent (2)\qquad Let $G\subseteq \bP^{\bar{g},\zeta}_f$ be generic over
$\bV$. For $\eta\in\fseo$ define
\[\begin{array}{l}
g^{\ell,*}_{\eta}(\zeta)=r^p_{\ell,\eta}\qquad\mbox{ where } p\in G\cap
\cI^{\lh(\eta)+1},\\
g^{\ell,*}_\eta(\xi)=g^\ell_\eta(\xi)\qquad\mbox{for }\xi<\zeta.
  \end{array}\]
It follows immediately from \ref{cl1} (and the definition of the order on
$\bP^{\bar{g},\zeta}_f$) that the above conditions define a
$\zeta+1$--approximation $\bar{g}^*=\langle g^{\ell,*}_\eta:\ell<2,\ \eta\in
{}^{\textstyle\omega{>}}\omega\rangle$ which agrees with $f$ and extends
$\bar{g}$. 
\end{proof}

\begin{theorem}
\label{getDa}
Assume that $\kappa$ is an uncountable cardinal such that $\kappa^{<\kappa}
=\kappa$. Then there is a ccc forcing notion $\bP$ of size $\kappa$ such
that
\[\forces_{\bP}\mbox{`` }(\circledast^{\rm Da})\ +\ \con=\kappa\ +\
\mbox{there is a spread family of size $\con$ ''.}\]
\end{theorem}

\begin{proof}
Using standard bookkeeping argument build inductively a finite support
iteration $\langle\bP_\alpha,\dbQ_\alpha:\alpha<\kappa\rangle$ and sequences
$\langle\zeta_\alpha:\alpha<\kappa\rangle$, $\langle\dot{\bar{g}}_\alpha:
\alpha<\kappa\rangle$ and $\langle \dot{f}_\alpha:\alpha<\kappa\rangle$ such
that:
\begin{enumerate}
\item $\bQ_0$ is the forcing notion adding $\kappa$ many Cohen reals
$\langle\rho_\xi:\xi<\kappa\rangle\subseteq\fseo$ (by finite approximations;
so, in $\bV^{\bQ_0}$, $\con=\kappa$ and the family
$\cF=\{\rho_\xi:\xi<\kappa\}$ is spread; we use it in the clauses below),
\item $\zeta_\alpha<\kappa$, $\dot{f}_\alpha$ is a $\bP_\alpha$--name for a
function from $(\zeta_\alpha+1)\times (\zeta_\alpha+1)$ to $\mbR$,
$\dot{\bar{g}}_\alpha$ is a $\bP_\alpha$--name for a
$\zeta_\alpha$--approximation (for the family $\cF$ added by $\bQ_0$) which
agrees with $\dot{f}_\alpha\rest (\zeta_\alpha\times\zeta_\alpha)$, 
\item $\forces_{\bP_{1+\alpha}}\dbQ_{1+\alpha}=\bP^{\dot{\bar{g}}_\alpha,
\zeta_\alpha}_{\dot{f}_\alpha}$ \qquad (for $\cF$),
\item if $\dot{f}$ is a $\bP_\kappa$ name for a function from
$(\zeta+1)\times (\zeta+1)$ to $\mbR$, $\zeta<\kappa$ and $\dot{\bar{g}}$ is
a $\bP_\kappa$--name for a $\zeta$--approximation which agrees with $\dot{f}
\rest(\zeta\times\zeta)$ then for some $\alpha<\kappa$, $\alpha>\omega$ we
have
\[\dot{\bar{g}}=\dot{\bar{g}}_\alpha,\quad \dot{f}=\dot{f}_\alpha,\quad
\zeta=\zeta_\alpha.\]
\end{enumerate}
Clearly $\bP_\kappa$ is a ccc forcing notion (with a dense subset) of size
$\kappa$. It follows from \ref{forissweet}(2), \ref{clear} that
$\forces_{\bP_\kappa}(\circledast^{\rm Da})$ (and clearly
$\forces_{\bP_\kappa}\con=\kappa$). Moreover, by \ref{forissweet}(1),
\ref{pressweet} we know that, in $\bV^{\bQ_0}$, for each $\alpha\in
[1,\kappa]$ the forcing notion $\bP_\alpha\rest [1,\kappa)$ is $\cF$--sweet,
so  
\[\forces_{\bP_\alpha}\mbox{`` $\cF$ is a spread family of size $\kappa$
''}\] 
(by \ref{swspr}). 
\end{proof}

\section{When $(\circledast^{\rm Da})$ fails.}
In this section we will strengthen the result of \cite[3.6]{Sh:675}
mentioned in \ref{675}(2) giving its combinatorial heart. 

\begin{definition}
\label{addinv}
\begin{enumerate}
\item For a function $h$ such that $\dom(h)\subseteq\cX\times\cY$ and
$\rng(h)\subseteq\cZ$ and a positive integer $n$ we define 
\[\begin{array}{ll}
\kappa(h,n)=\min\{|\cA_0|+|\cA_1|:& \cA_0\subseteq\cP(\cX)\ \&\ \cA_1
\subseteq\cP(\cY)\ \&\\
&\quad (\forall w\in [\cX]^{\textstyle n})(\exists A\in\cA_0)(w\subseteq A)\
\&\\
&\quad (\forall w\in [\cY]^{\textstyle n})(\exists A\in\cA_1)(w\subseteq A)\
\&\\  
&\quad (\forall A_0\in\cA_0)(\forall A_1\in \cA_1)(h[A_0\times A_1]\neq\cZ)
\}.
\end{array}\]
If $\cX=\cY$ and $h $ is as above, and $n$ is a positive integer then we
define  
\[\begin{array}{lr}
\kappa^-(h,n)=\min\{|\cA|:& \cA\subseteq\cP(\cX)\ \&\ (\forall w\in [\cX]^{
\textstyle n})(\exists A\in\cA)(w\subseteq A)\ \&\quad\\
&\quad (\forall A\in\cA)(h[A\times A]\neq\cZ)\ \}.
\end{array}\]
\item For $\bar{c}=\langle c_n:n<\omega\rangle\in{}^{\textstyle\omega}\mbR$
and $\bar{d}=\langle d_n:n<\omega\rangle\in{}^{\textstyle\omega}\mbR$ let 
$h^\oplus(\bar{c},\bar{d})=\sum\limits_{n<\omega} c_n\cdot d_n$ (defined if
the series converges). 
\end{enumerate}
\end{definition}

We will deal with the following variant of the property $(\circledast^{\rm
Da})$.  

\begin{definition}
\label{moreDa}
For a function $h:\sqR\times\sqR\longrightarrow\mbR$ let $(\circledast^{\rm
Da}_h)$ mean: 
\begin{enumerate}
\item[$(\circledast^{\rm Da}_h)$] For each $f:\mbR\times\mbR\longrightarrow
\mbR$ there are functions $g^0_n,g^1_n:\mbR\longrightarrow\mbR$ (for $n<
\omega$) such that 
\[(\forall x,y\in\mbR)\big(f(x,y)=h(\langle g^0_n(x):n<\omega\rangle,\
\langle g^1_n(y):n<\omega\rangle)\big).\] 
\end{enumerate}
(So $(\circledast^{\rm Da})$ is $(\circledast^{\rm Da}_{h^\otimes})$, where
$h^\oplus$ is as defined in \ref{addinv}(2).)
\end{definition}

\begin{proposition}
\label{notDaone}
Assume that a function $h:\sqR\times\sqR\longrightarrow\mbR$ is such that
on of the following condition holds:
\begin{enumerate}
\item[(A)] $\kappa(h,1)<2^{\kappa(h,1)}=\con$,\qquad or
\item[(B)] $\kappa(h,1)\leq\mu<\con$ for some regular cardinal $\mu$,\qquad
or
\item[(C)] $\kappa^-(h,2)\leq\mu<\con$ for some regular cardinal $\mu$.
\end{enumerate}
Then $(\circledast^{\rm Da}_h)$ fails.
\end{proposition}

\begin{proof}
First let us consider the case of the assumption (A).
Let $\cA_0,\cA_1\subseteq\cP(\sqR)$ exemplify the minimum in the definition
of $\kappa(h,1)$, $\cA_\ell=\{A^\ell_\xi:\xi<\kappa(h,1)\}$ (we allow
repetitions). Choose a sequence $\langle r_\xi:\xi<\kappa(h,1)\rangle$ of
pairwise distinct reals and fix enumerations $\langle s_\varepsilon:
\varepsilon<\con\rangle$ of $\mbR$ and $\langle\varphi_\varepsilon:
\varepsilon<\con\rangle$ of ${}^{\textstyle\kappa(h,1)}\kappa(h,1)$. Let
$f:\mbR\times\mbR\longrightarrow\mbR$ be such that
\[(\forall\xi<\con)(\forall\xi<\kappa(h,1))\big(f(s_\varepsilon,r_\xi)
\notin h[A^0_\xi\times A^1_{\varphi_\varepsilon(\xi)}]\big).\]
We claim that the function $f$ witnesses the failure of $(\circledast^{\rm
Da}_h)$. So suppose that $g^0_n,g^1_n:\mbR\longrightarrow\mbR$. For
$\xi<\kappa(h,1)$ let $\bar{b}_\xi=\langle g^1_n:n<\omega\rangle\in\sqR$ and
let $\varphi(\xi)<\kappa(h,1)$ be such $\bar{b}_\xi\in A^1_{\varphi(\xi)}$. 
Take $\varepsilon<\con$ such that $\varphi=\varphi_\varepsilon$ and let
$\bar{a}_\varepsilon=\langle g^0_n(\varepsilon):n<\omega\rangle$. Fix $\xi^*
<\kappa(h,1)$ such that $\bar{a}_\varepsilon\in A^0_{\xi^*}$ and note that
$h(\bar{a}_\varepsilon,\bar{b}_{\xi^*})\in h[A^0_{\xi^*}\times
A^1_{\varphi_\varepsilon(\xi_0)}]$, so 
\[f(s_\varepsilon,r_{\xi^*})\neq h(\bar{a}_\varepsilon,\bar{b}_{\xi^*})=
h(\langle g^0_n(s_{\varepsilon}):n<\omega\rangle,\ \langle g^1_n(r_{\xi^*}):
n<\omega\rangle).\]
\medskip

\noindent Suppose now that we are in the situation (B). Let $c_0,c_1:\mu^+
\times\mu^+\longrightarrow\kappa(h,1)$ be such that for any sets $X_0,X_1
\in[\mu^+]^{\textstyle\mu^+}$ we have
\[(\forall\zeta_0,\zeta_1<\kappa(h,1))(\exists\langle\varepsilon_0,
\varepsilon_1\rangle\in X_0\times X_1)(c_0(\varepsilon_0,\varepsilon_1)=
\zeta_0\ \&\ c_1(\varepsilon_0,\varepsilon_1)=\zeta_1)\]
(see e.g.\ \cite[ch III]{Sh:g}). Let $\cA_0,\cA_1\subseteq\cP(\sqR)$
exemplify $\kappa(h,1)$, $\cA_\ell=\{A^\ell_\zeta:\zeta<\kappa(h,1)\}$
(with possible repetitions). Choose a sequence $\langle r_\varepsilon:
\varepsilon<\mu^+\rangle$ of pairwise distinct reals and a function
$f:\mbR\times\mbR\longrightarrow\mbR$ such that
\[(\forall\varepsilon_0,\varepsilon_1<\mu^+)\big(f(r_{\varepsilon_0},
r_{\varepsilon_1})\notin h[A^0_{c_0(\varepsilon_0,\varepsilon_1)}\times
A^1_{c_1(\varepsilon_0,\varepsilon_1)}]\big).\]
Now suppose that $g^0_n,g^1_n:\mbR\longrightarrow\mbR$ and let
$\bar{a}^\ell_\varepsilon=\langle g^\ell_n(r_\varepsilon):n<\omega\rangle$. 
Choose $X_0,X_1\in [\mu^+]^{\textstyle\mu^+}$ and $\zeta_0,\zeta_1<\kappa(h,
1)$ such that $\bar{a}^\ell_\varepsilon\in A^\ell_{\zeta_\ell}$ whenever
$\varepsilon\in X_\ell$. Take $\varepsilon_\ell\in X_\ell$ (for $\ell<2$)
such that $c_0(\varepsilon_0,\varepsilon_1)=\zeta_0$, $c_1(\varepsilon_0,
\varepsilon_1)=\zeta_1$. Then $h(\bar{a}^0_{\varepsilon_0},\bar{a}^1_{
\varepsilon_1})\in h[A^0_{c_0(\varepsilon_0,\varepsilon_1)}\times A^1_{c_2
(\varepsilon_0,\varepsilon_1)}]$, so $f(r_{\varepsilon_0},
r_{\varepsilon_1}) \neq h(\langle g^0_n(r_{\varepsilon_0}):n<\omega\rangle,\
\langle g^1_n(r_{\varepsilon_1}):n<\omega\rangle)$. 
\medskip

\noindent Now, suppose that the assumption (C) holds. Let $\{A_\xi:\xi<
\kappa^-(h,2)\}$ be a family witnessing the minimum in the definition of
$\kappa^-(h,2)$. Take a function $c:\mu^+\times\mu^+\longrightarrow\kappa^-
(h,2)$ such that for every $X\in [\mu^+]^{\textstyle\mu^+}$ and $\zeta<
\kappa^-(h,2)$ there are $\varepsilon_0<\varepsilon_1$, both in $X$, such
that $c(\varepsilon_0,\varepsilon_1)=\zeta$ (see e.g.\ \cite[ch III]{Sh:g}).
Take a sequence $\langle r_\varepsilon:\varepsilon<\mu^+\rangle$ of distinct
reals and define a function $f:\mbR\times\mbR\longrightarrow\mbR$ so that
\[(\forall\varepsilon_0,\varepsilon_1<\mu^+)(f(r_{\varepsilon_0},r_{
\varepsilon_1})\notin h [A_{c(\varepsilon_0,\varepsilon_1)}\times A_{c(
\varepsilon_0,\varepsilon_1)}]).\]
Like before, suppose that $g^0_n,g^1_n:\mbR\longrightarrow\mbR$ and let
$\bar{a}^\ell_\varepsilon=\langle g^\ell_n(r_\varepsilon):n<\omega\rangle$. 
For each $\varepsilon<\mu^+$ there is $\zeta_\varepsilon\in\kappa^-(h,2)$
such that $\bar{a}^0_\varepsilon,\bar{a}^1_\varepsilon\in A_{\zeta_{
\varepsilon}}$. Take a set $X\in [\mu^+]^{\textstyle \mu^+}$ and $\zeta^*<
\kappa^-(h,2)$ such that $(\forall\varepsilon\in X)(\zeta_\varepsilon=
\zeta^*)$. Then choose $\varepsilon_0<\varepsilon_1$ both in $X$ so that
$c(\varepsilon_0,\varepsilon_1)=\zeta^*$. By our choices, $\bar{a}^0_{
\varepsilon_0},\bar{a}^1_{\varepsilon_1}\in A_{c(\varepsilon_0,
\varepsilon_1)}$ and $h(\bar{a}^0_{\varepsilon_0},\bar{a}^1_{\varepsilon_1})
\in A_{c(\varepsilon_0,\varepsilon_1)}$. But this implies that $h(\langle
g^0_n(r_{\varepsilon_0}):n<\omega\rangle,\ \langle g^1_n(r_{\varepsilon_1}):
n<\omega\rangle)\neq f(r_{\varepsilon_0},r_{\varepsilon_1})$.
\end{proof}

Now the phenomenon of \cite[3.6]{Sh:675} is described in a combinatorial way
by \ref{notDaone}, if one notices the following observation.

\begin{proposition}
Let $h:\sqR\times\sqR\longrightarrow\sqR$ be a function with an absolute
definition (with parameters from the ground model). Suppose that $\bP=
\langle\bP_\alpha,\dbQ_\alpha:\alpha<\omega_1\rangle$ is a finite support
iteration of non-trivial forcing notions. Then for each $0<n<\omega$
\[\forces_{\bP_{\omega_1}}\kappa(h,n)=\kappa^-(h,n)=\aleph_1.\]
\end{proposition}

\begin{proof}
Work in $\bV^{\bP_{\omega_2}}$. For $\alpha<\omega_1$ let $A_\alpha=\bV^{
\bP_\alpha}\cap \sqR$. Clearly $\sqR=\bigcup\limits_{\alpha<\omega_1}
A_\alpha$ and for each $\alpha,\beta<\omega_1$ we have $h[A_\alpha\times
A_\beta]\neq\sqR$ (remember that the function $h$ has definition with
parameters in the ground model; at each limit stage of the iteration Cohen
reals are added).  
\end{proof}

\section{Concluding remarks}
One can notice some similarities between the property $(\circledast)^{\rm
Da}$ and the rectangle problem.
\begin{definition}
\begin{enumerate}
\item Let $\cR_2$ be the family of all rectangles in $\mbR\times\mbR$, i.e.\
sets of the form $A\times B$ for some $A,B\subseteq\mbR$. Let $\cB(\cR_2)$
be the $\sigma$--algebra of subsets of $\mbR\times\mbR$ generated by the
family $\cR_2$ and let $\cB_\alpha(\cR_2)$ be defined inductively by:
$\cB_0(\cR_2)$ consists of all elements of $\cR_2$ and their complements,
$\cB_\alpha(\cR_2)=\bigcup\limits_{\beta<\alpha}\cB_\beta(\cR_2)$ for limit
$\alpha$, and $\cB_{\alpha+1}(\cR_2)$ is the collection of all countable
unions $\bigcup\limits_{n<\omega} A_n$ such that each $A_n$ is in
$\cB_\alpha(\cR_2)$ and of the complements of such unions. (So $\cB(\cR_2)=
\cB_{\omega_1}(\cR_2)$.)
\item Let us introduce the following properties of the family of subsets of
$\mbR\times\mbR$: 
\begin{enumerate}
\item[$(\boxdot^{\rm Ku})$]\qquad $\cP(\mbR\times\mbR)=\cB(\cR_2)$, 
\item[$(\boxdot^{\rm Ku}_\alpha)$]\qquad $\cP(\mbR\times\mbR)=\cB_\alpha(
\cR_2)$
\end{enumerate}
\end{enumerate}
\end{definition}

Kunen \cite[\S 12]{Ku68} showed the following.

\begin{theorem}
\begin{enumerate}
\item (See \cite[Thm 12.5]{Ku68})\quad ${\bf MA}$ implies $(\boxdot^{\rm
Ku}_2)$.  
\item (See \cite[Thm 12.7]{Ku68})\quad  If $\bP$ is the forcing notion for
adding $\aleph_2$ Cohen reals then $\forces_{\bP}\neg (\boxdot^{\rm Ku})$.
\end{enumerate}
\end{theorem}

The relation between $(\circledast^{\rm Da})$ and $(\boxdot^{\rm Ku})$ is
still unclear, though the first implies the second.

\begin{proposition}
\qquad$(\circledast^{\rm Da})\quad\Rightarrow\quad(\boxdot^{\rm Ku}_\omega)$
\end{proposition}

\begin{proof}
Suppose that $A\subseteq\mbR\times\mbR$ and let $f:\mbR\times\mbR
\longrightarrow 2$ be it characteristic function. Let $g^0_n,g^1_n$ be given
by $(\circledast^{\rm Da})$ for the function $f$. For a rational number $q$,
$n<\omega$ and $\ell<2$ put 
\[A^\ell_{q,n}\stackrel{\rm def}{=}\{x\in\mbR: g^\ell_n(x)<q\}.\]
It should be clear that the set $A$ can be represented as a Boolean
combination of finite depth of rectangles $A^0_{q,n}\times A^1_{q',n}$ (we
do not try to safe on counting the quantifiers).
\end{proof}

The following questions arise naturally in this context.

\begin{problem}
\begin{enumerate}
\item Does $(\boxdot^{\rm Ku}_\omega)$ (or $(\boxdot^{\rm Ku})$) imply
$(\circledast^{\rm Da})$?
\item Is it consistent that for some countable limit ordinal $\alpha$ we
have $(\boxdot^{\rm Ku}_{\alpha+1})$ but $(\boxdot^{\rm Ku}_{\alpha})$
fails? 
\end{enumerate}
\end{problem}


\begin{thebibliography}{BaJu95}
\makeatletter \renewcommand{\@biblabel}[1]{[#1]} \makeatother

\bibitem[BaJu95]{BaJu95}Tomek Bartoszy\'nski and Haim Judah.
\newblock {\em {Set Theory: On the Structure of the Real Line}}.
\newblock A K Peters, Wellesley, Massachusetts, 1995.

\bibitem[Ci97]{Ci97}K.~Ciesielski.
\newblock Set theoretic real analysis.
\newblock {\em J. of Applied Analysis}, {\bf 3}:143--190, 1997.

\bibitem[Da74]{Da74}R.~O. Davies.
\newblock {Representation of functions of two variables as sums of rectangular
  functions I}.
\newblock {\em Fundamenta Mathematicae}, pages 177--183, 1974.

\bibitem[Ku68]{Ku68}Kenneth Kunen.
\newblock {\em {Inaccessibility properties of cardinals}}.
\newblock PhD thesis, Stanford University, 1968.

\bibitem[Mi91]{Mi91}Arnold~W. Miller.
\newblock {Arnie Miller's problem list}.
\newblock In Haim Judah, editor, {\em {Set Theory of the Reals}}, volume~6 of
  {\em Israel Mathematical Conference Proceedings}, pages 645--654.
\newblock {Proceedings of the Winter Institute held at Bar--Ilan University,
  Ramat Gan, January 1991}.

\bibitem[Mixx]{Mixx}Arnold~W. Miller.
\newblock {Some interesting problems}.
\newblock {circulated notes; available at {\tt
  http://www.math.wisc.edu/$\sim$miller}}.

\bibitem[Sh:98]{Sh:98}Saharon Shelah.
\newblock {Whitehead groups may not be free, even assuming CH. II}.
\newblock {\em {Israel Journal of Mathematics}}, {\bf 35}:257--285, 1980.

\bibitem[Sh:g]{Sh:g}Saharon Shelah.
\newblock {\em {Cardinal Arithmetic}}, volume~29 of {\em {Oxford Logic
  Guides}}.
\newblock {Oxford University Press}, 1994.

\bibitem[Sh 675]{Sh:675}Saharon Shelah.
\newblock {On Ciesielski's Problems}.
\newblock {\em {Journal of Applied Analysis}}, {\bf 3}(2):{191--209}, 1997.

\end{thebibliography}
\def\germ{\frak} \def\scr{\cal}
  \ifx\documentclass\undefinedcs\def\rm{\fam0\tenrm}\fi
  \def\defaultdefine#1#2{\expandafter\ifx\csname#1\endcsname\relax
  \expandafter\def\csname#1\endcsname{#2}\fi} \defaultdefine{Bbb}{\bf}
  \defaultdefine{frak}{\bf} \defaultdefine{mathbb}{\bf}
  \defaultdefine{beth}{BETH} \def\bbfI{{\Bbb I}} \def\mbox{\hbox}
  \def\text{\hbox} \def\om{\omega} \def\Cal#1{{\bf #1}} \def\pcf{pcf}
  \def\restriction{{|}} \def\club{CLUB} \def\w{\omega} \def\exist{\exists}
  \def\se{{\germ se}} \def\bb{{\bf b}} \def\equivalence{\equiv}
  \def\cite#1{[#1]}

\end{document}